\documentclass[leqno,a4paper,oneside]{amsart}
\usepackage{amsmath,amssymb,amsfonts,latexsym,amscd,mathrsfs,mathtools}
\usepackage{graphicx,comment,multirow,color,caption}

\usepackage{mathptmx,courier}
\usepackage[scaled=1.0]{helvet}

\usepackage[
  colorlinks,
  citecolor=blue,
  linkcolor=blue,
  urlcolor=blue
]{hyperref}

\usepackage[
  width= 135mm,
]{geometry}

\usepackage{wrapfig}
\usepackage[%
    font={small,sf},
    labelfont=normalfont,
    format=hang,    
    format=plain,
    margin=0pt,
    width=0.8\textwidth,
]{caption}
\usepackage[list=true]{subcaption}

\title[Simultaneous crossing number]{Multiple points of view: 
The simultaneous crossing number \\ for knots with doubly transvergent diagrams}

\author{Christoph Lamm}
\address{R\"{u}ckertstr. 3, 65187 Wiesbaden, Germany}
\email{christoph.lamm@web.de}

\author{Michael Eisermann}
\address{Fachbereich Mathematik, Universität Stuttgart, Germany}
\email{Michael.Eisermann@mathematik.uni-stuttgart.de}

\theoremstyle{plain}
  \newtheorem{theorem}{Theorem}[section]

\theoremstyle{definition}
  \newtheorem{definition}[theorem]{Definition}
  \newtheorem{question}[theorem]{Question}
  \newtheorem{project}[theorem]{Project}
  \newtheorem{remark}[theorem]{Remark}
  \newtheorem{example}[theorem]{Example}
  \newtheorem{conjecture}[theorem]{Conjecture}

\renewcommand{\S}{\mathbb{S}}
\newcommand{\cn}{\operatorname{cr}}
\newcommand{\scn}{\operatorname{sim}}
\newcommand{\isoto}{\overset{\!\sim}{\to}}

\begin{document} 

\begin{abstract}
The simultaneous crossing number is a new knot invariant which is defined for strongly 
invertible knots having diagrams with two orthogonal transvergent axes of strong inversions.
Because the composition of the two inversions gives a cyclic period of order 2 with an axis 
orthogonal to the two axes of strong inversion, knot diagrams with this property have three 
characteristic orthogonal directions. 
We define the \textit{simultaneous crossing number}, $\scn(K)$, as the minimum of the 
sum of the numbers of crossings of projections in the 3 directions, where the
minimum is taken over all embeddings of $K$ satisfying the symmetry condition.
Dividing the simultaneous crossing number by the usual crossing number, $\cn(K)$, of a knot gives 
a number $\ge 3$, because each of the 3 diagrams is a knot diagram of the knot in question.
We show that $\liminf_{\cn(K) \to \infty} \scn(K)/\cn(K) \le 8$, when the minimum over all knots
and the limit over increasing crossing numbers is considered.
\end{abstract}



\captionsetup{belowskip=8pt,aboveskip=8pt}

\reversemarginpar


\maketitle

\section{Introduction} \label{sec:Introduction}

Let $K \subset \S^3$ be a (polygonal or smooth) knot. 
A \emph{strong inversion} $\rho$ of the knot $K$ 
is a smooth involution $\rho \colon (\S^3,K) \isoto (\S^3,K)$ preserving the orientation on $\S^3$ while reversing it on $K$.

\begin{figure}[hbtp]
\centering
\includegraphics[scale=0.8]{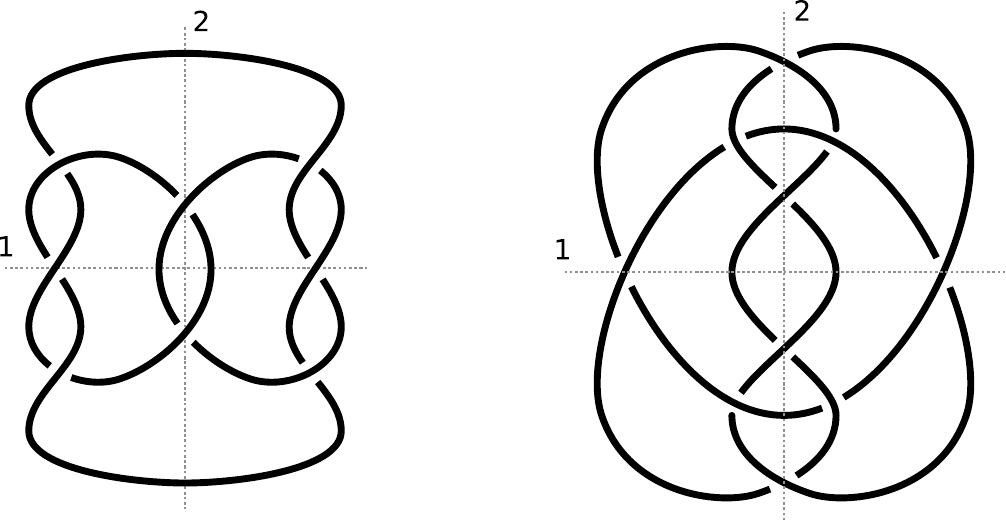}
\caption{Examples of doubly transvergent diagrams for the knots $8_5$ and $10_{122}$}
\label{doubly_transvergent_diagrams}
\end{figure}

Consider, for example, a knot $K$ as in the diagrams of Figure \ref{doubly_transvergent_diagrams}.
The half rotation around the $x$-axis is a strong inversion $\rho_x$,
likewise $\rho_y$ around the $y$-axis, orthogonal to the $x$-axis.
We call such a knot diagram \textit{doubly transvergent}.

The composition $\rho_z = \rho_x \rho_y = \rho_y \rho_x$ 
is a half rotation around the $z$-axis, perpendicular to the diagram plane, 
mapping $K$ to itself and preserving its orientation.
This shows that every knot $K$ with a doubly transvergent diagram 
automatically has cyclic period 2.

Knots with doubly transvergent diagrams therefore 
have three characteristic orthogonal directions.
Projecting along the two axes of strong inversion,
onto the $y$-$z$-plane or the $x$-$z$-plane,
results in two intravergent diagrams of the knot (this terminology was introduced in \cite{Boyle}).
In this articles we study the sum of the crossing numbers 
of the three diagrams obtained in this way.

\begin{definition}
Let $K \subset \S^3$ be a knot with a doubly transvergent diagram.

The \textit{simultaneous crossing number}, denoted $\scn(K)$, 
is the minimum of the sum of the three crossing numbers,
obtained by projecting along the three characteristic axes.

Here the minimum  is taken over all embeddings of $K$ 
giving a doubly transvergent diagram in one direction
and two intravergent diagrams in the other two directions.
\end{definition}

We have $\scn(K) \ge 3 \cdot \cn(K)$, where $\cn(K)$ is the crossing number of $K$.

\begin{question} \label{quest:Parity}
Is $\scn(K) = 3 \cdot \cn(K)$ possible?
\end{question}

We note the following parity properties of the numbers of crossings 
in transvergent and intravergent diagrams:

\begin{remark}
A (simply) transvergent diagram 
may have an even or an odd number of crossings.
For every doubly transvergent diagram, however,
the crossing number is always even, because 
the rotation $\rho_z$ pairs crossings.
Therefore, a doubly transvergent diagram cannot realize 
the knot's minimal crossing number, if $\cn(K)$ is odd.
\end{remark}

\begin{remark}
In the other two directions of projection, we have intravergent diagrams. Apart from the central crossing, all other crossings are paired 
by the rotation around the center, so the crossing number is always odd.
Therefore, an intravergent diagram cannot realize the knot's minimal crossing number, if $\cn(K)$ is even.
\end{remark}

Summing over the three directions, we answer Question \ref{quest:Parity} negatively:
\begin{align}
\label{eq:LowerBound}
\scn(K) \ge \begin{cases}
3 \cdot \cn(K) + 1, & \text{if $\cn(K)$ is odd,} \\ 
3 \cdot \cn(K) + 2, & \text{if $\cn(K)$ is even.}
\end{cases}
\end{align}

As an easy example, we calculate the simultaneous crossing number of the trivial knot.
This is currently the only knot whose simultaneous crossing number is known.

\begin{example} \label{ex:TrivialKnot}
The simultaneous crossing number of the trivial knot $T$ is $\scn(T) = 2$.
\end{example}

\begin{figure}[hbtp]
\centering
\includegraphics[scale=1.6]{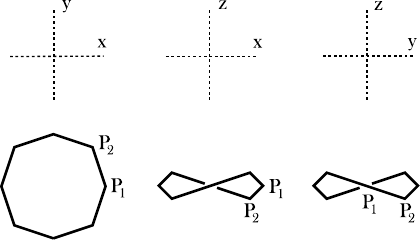}
\caption{The three views of a doubly transvergent representation of the trivial knot}
\label{trivial_knot}
\end{figure}

\begin{proof}
Since $\cn(T) = 0$, the lower bound 
\eqref{eq:LowerBound} yields $\scn(T) \ge 2$.
In the other direction, the closed polygon with vertices 
$P_1(4,0,0)$, $P_2(3,3,-1)$, $P_3(0,4,0)$, $P_4(-3,3,1)$,
$P_5(-4,0,0)$, $P_6(-3,-3,-1)$, $P_7(0,-4,0)$ and $P_8(3,-3,1)$ 
satisfies the symmetry  conditions and has a total number 
of 2 crossings, as shown in Figure \ref{trivial_knot}.
\end{proof}

Besides determining the simultaneous crossing number for individual knots, there are also other
topics to consider. Here is the plan of the article.

\subsection*{Plan of the article.}
In section \ref{asymptotic_behaviour} we ask for the lower and upper values of the quotients
$\scn(K)/\cn(K)$ for increasing $\cn(K)$. We show that for twist knots $\scn(K) \le 8 \cdot \cn(K) - 10$
and therefore find $\liminf_{\cn(K) \to \infty} \scn(K)/\cn(K) \le 8$ when all knot $K$ with doubly
transvergent diagrams are considered. (For the limit superior we do not yet know any bound.)

In section \ref{related_concepts} we discuss concepts which are related to the simultaneous crossing number,
for instance the super crossing index.

In section \ref{equivalence_classes} we investigate which types of doubly transvergent diagrams
occur with respect to the two equivalence classes of strong inversions, and formulate a conjecture.

The final section \ref{summary} contains a list of more detailed questions and conjectures.

We add three appendices: One on the symmetries of chord diagrams for doubly transvergent diagrams, the second one on using the software `Blender' for our purpose, and the final one with 
supplements on the remaining twist knot cases, not treated in section \ref{asymptotic_behaviour}.

\section{The limit behaviour for increasing crossing number} \label{asymptotic_behaviour}

The simultaneous crossing number is currently known 
only for the trivial knot (Example \ref{ex:TrivialKnot}).
For the trefoil knot $3_1$ we have $\cn(3_1) = 3$,
so the lower bound \eqref{eq:LowerBound} yields $\scn(3_1) \ge 10$.

\begin{figure}[hbtp]
\includegraphics[scale=0.2]{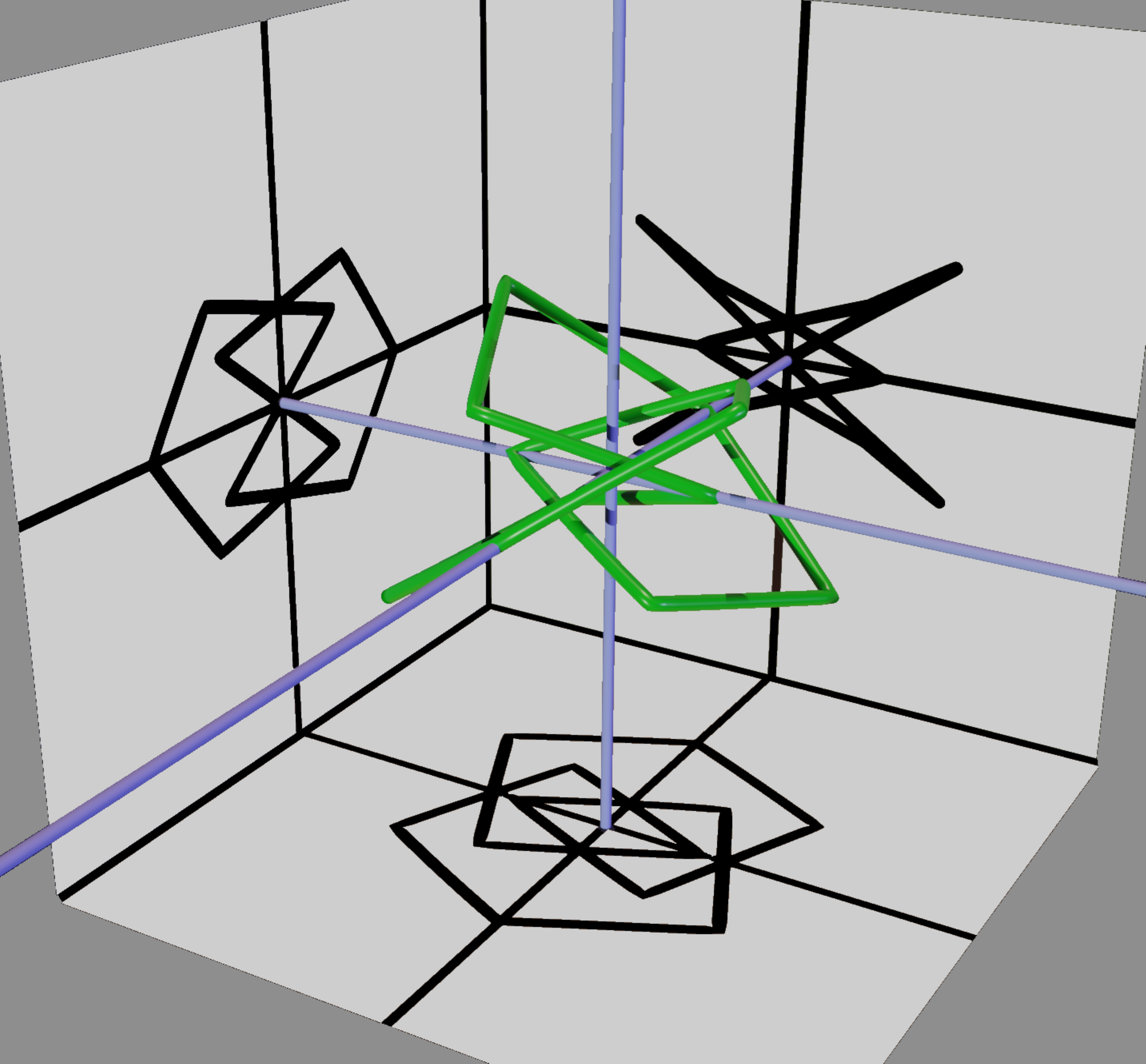}
\caption{A trefoil with a doubly transvergent diagram in the $x$-$y$-plane and intravergent diagrams in the other two coordinate planes}
\label{trefoil_with_axes}
\end{figure}

Experimentally, we find the upper bound $\scn(3_1) \le 14$: 
The closed polygon in Figure  \ref{trefoil_with_axes} 
has $4$ crossings in the transvergent diagram (projecting along the $z$-axis), plus $3$ resp.\ $7$ crossings in the intravergent diagrams (projecting along the $x$-axis and the $y$-axis).

We do not give detailed coordinates for this polygonal representation, 
because the trefoil is part of the twist knot family, 
which is described in detail in the following paragraphs.

\begin{question} \label{trefoil_sim_c_question}
Is $\scn(3_1)=14$?
\end{question}

We cannot yet answer individual questions as this one.
A possible \emph{brute force} attack would be to systematically 
enumerate all cases with a total number of crossings less than $14$, say.
While this seems doable, we currently do not know how to carry this out.

Instead, we would like to study the lower and upper values of the quotients $\scn(K)/\cn(K)$ for knot families with increasing $\cn(K)$.
In the example of twist knots, we observe, that for every integer $n \ge 3$ there is a knot with a doubly symmetric diagram.
Therefore, we can define $\scn_{\min}(n)$ and $\scn_{\max}(n)$ as the smallest (largest) simultaneous crossing number of all knots with doubly symmetric 
diagrams $K$ with $\cn(K)=n$, for $n \ge 3$. We know, that $\scn_{\min}(n) \ge 3 \cdot n$, and hence $\liminf_{n \to \infty} \scn_{\min}(n)/n \ge 3$.

\begin{question} \label{liminf_question}
What is the value of $\liminf_{n \to \infty} \scn_{\min}(n)/n$?
\end{question} 

The analogous question for the maximal values is then:

\begin{question} \label{limsup_question}
Is $\limsup_{n \to \infty} \scn_{\max}(n)/n < \infty$? If yes, what is its value?
\end{question} 

For the example family of twist knots ($3_1$, $4_1$, $5_2$, $6_1$, $\ldots$) we find an upper bound for 
the simultaneous crossing number and prove the following theorem:

\begin{theorem} \label{twist_knot_theorem}
Let $K_n$ be the twist knot with crossing number $n$.
\newline
Then $\scn(K_n) \le 8 \cdot n-10$ and therefore $\liminf_{n \to \infty} \scn_{\min}(n)/n \le 8$.
\end{theorem} 

\begin{figure}[hbtp]
\centering
\includegraphics[scale=1.3]{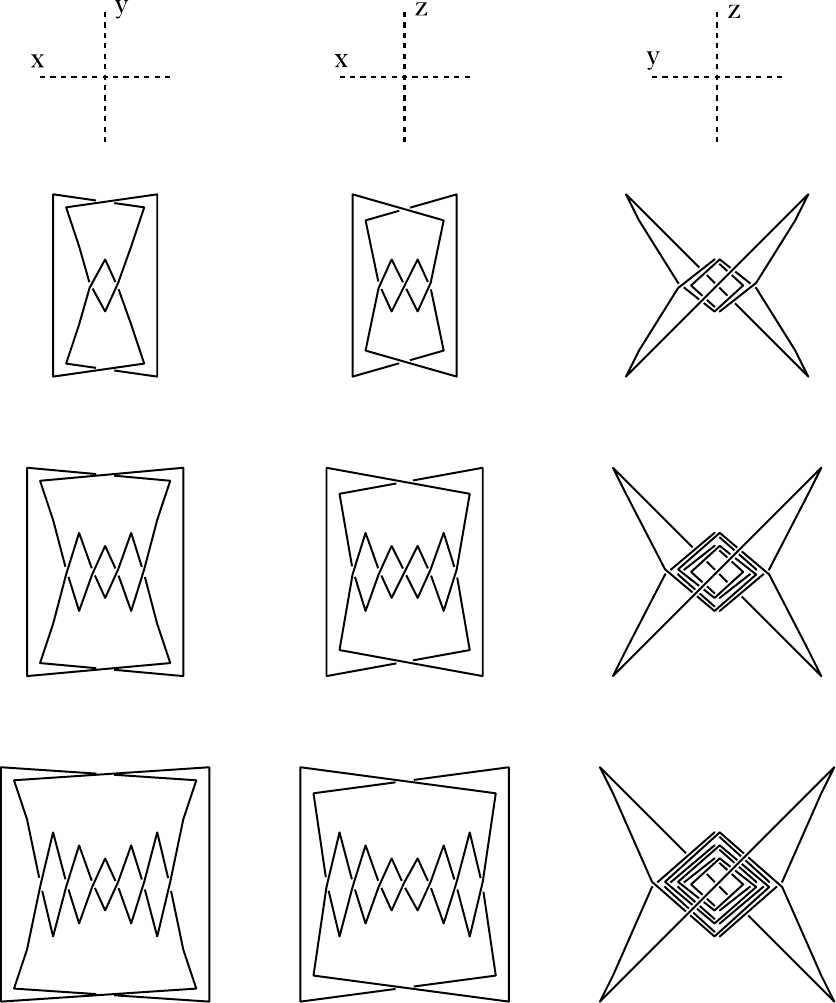}
\caption{The twist knots $C(2,2)$, $C(4,2)$ and $C(6,2)$ in three orthogonal projections}
\label{projections_twist_knots_positive}
\end{figure}

\begin{proof}
Using the Conway notation $C(a_1, \ldots, a_r)$ for two-bridge knots, we obtain doubly transvergent
diagrams for twist knots: a) The family with even crossing numbers is described by $C(2k,2)$ for $k \ge 1$.
This family consists of the knots $4_1$, $6_1$, $8_1$, etc. b) The knots $C(2k,-2)$ have odd crossing
numbers, and we get $3_1$, $5_2$, $7_2$, etc.

We give explicit polygonal versions for these knots. They are shown in Figures \ref{projections_twist_knots_positive}
and \ref{projections_twist_knots_negative} for $k=1,2,3$. The structures of the polygons are quite similar, but
there are differences between $C(2k,2)$ and $C(2k,-2)$, and also between even and odd parameters $k$. Therefore, we describe only one of the four sub-families in detail: $C(2k,2)$ for odd $k$, yielding $4_1$, $8_1$, etc.
Since the description for the other three sub-families is similar, we relegate it to Appendix C.

It is enough to specify one fourth of the polygon, since the knot can be built from that part, using its symmetry.
As the smallest example, we start with $C(2,2)$. 

One fourth of the polygon for $C(2,2)$ is described as follows: $P_1(0,2,0)$, $P_2(1,-\delta,-2)$, 
$P_3(2,-3,\delta)$, $P_4(3,-6,5)$, $P_5(-4,-7,7)$, $P_6(-4,0,0)$. 
The parameter $\delta > 0$ ensures that near the axes there is no collision of vertices in the three orthogonal 
projections, when the complete polygon is constructed. We usually choose $\delta=0.2$.
After generating the other 3 polygonal parts, and gluing them together, the final polygon has 20 vertices 
($=24-4$, because 4 vertices are identified in the gluing with 4 other vertices).

Replacing $k$ by $k+2$, the spiral part around the $x$-axis is enlarged. 
The following algorithm in pseudo code gives the coordinates for $C(2k,2)$ in the general case for odd $k$. This is illustrated for $k=3$ in Figure \ref{polygon_one_fourth}.

\bigskip
\begin{figure}[hbtp]
\centering
\includegraphics[scale=1.3]{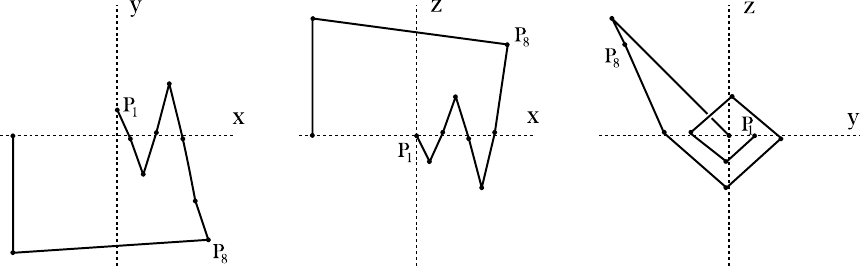}
\caption{Illustration of the polygonal curve for case 1. One fourth of the knot $C(6,2)$ consists of 10 vertices. After gluing, the knot has 36 vertices.}
\label{polygon_one_fourth}
\end{figure}

\footnotesize
\begin{verbatim}
Case 1: C(2k,2) for odd k

input m    # Examples: m=0 and k=1 give C(2,2)
k = 2*m+1  #           m=1 and k=3 give C(6,2)
d = 0.2    # The parameter delta
  
P[1] = (0, 2, 0)
P[2] = (1,-d,-2)

j = 2      # counter for vertices, starting with P[3]

for i = 0 to m-1:  # spiral part enlargement
    x = 4*i+2; y = -(2*i+3); z =  d; j = j+1; P[j] = (x,y,z)    		
    x = x+1;   z =  -y;      y =  d; j = j+1; P[j] = (x,y,z)		
    x = x+1;   y = z+1;      z = -d; j = j+1; P[j] = (x,y,z)		
    x = x+1;   z =  -y;      y = -d; j = j+1; P[j] = (x,y,z)
end

j = j+1; P[j] = ( 2*k  ,-k-2,  d)
j = j+1; P[j] = ( 2*k+1,-k-5,k+4)
j = j+1; P[j] = (-2*k-2,-k-6,k+6)
j = j+1; P[j] = (-2*k-2,   0,  0)
\end{verbatim}

\normalsize
\medskip
\noindent
Note, that all vertices contain the parameter $\delta$, with the exception of the first and the last three.
To finish the proof, we have to show, that
\begin{itemize}
	\item The definition of the fourth part yield knots with diagrams as shown in Figure \ref{projections_twist_knots_positive}.
	\item The sum of the numbers of crossings is $8 \cdot n-10$.
	\item Similar constructions also work for $C(2k,2)$ for even $k$, and for $C(2k,-2)$, and also yield $8 \cdot n-10$ crossings. This is relegated to Appendix C.
\end{itemize}

Some of the details will be omitted, but we are more specific for the calculation of the numbers of crossings:
In Table \ref{tab:numbers_crossings} we list the crossing numbers in the three directions (denoted by $p_z$, $p_y$ and $p_x$) and find that in each case $p_z+p_y=2 \cdot n+1$. In the $x$-direction the number is $p_x=6 \cdot n -11$, yielding a sum of $8 \cdot n-10$. 

The correctness of these values can be shown by induction, starting with the simplest cases and noting that replacing $k$ by $k+2$ results in increments of 4 crossings for the numbers $n$, $p_z$ and $p_y$, and of 24 crossings for $p_x$. 
\end{proof}

\begin{table}[hbtp]
  \centering
  \begin{tabular}{l|r|rrr|r}
  $K$ & $n$ & $p_z$  & $p_y$ & $p_x$ & sum \\
  \hline
  $C(2,-2)$ & 3 & 4 & 3 &  7 & 14\\
  $C(2,2)$  & 4 & 4 & 5 & 13 & 22\\
  $C(4,-2)$ & 5 & 6 & 5 & 19 & 30\\
  $C(4,2)$  & 6 & 6 & 7 & 25 & 38\\
  $C(6,-2)$ & 7 & 8 & 7 & 31 & 46\\
  $C(6,2)$  & 8 & 8 & 9 & 37 & 54\\
  \end{tabular}
  \caption{Table of the numbers of crossings in the twist knot family for $k \le 3$}
  \label{tab:numbers_crossings}
\end{table}

\begin{figure}[hbtp]
\centering
\includegraphics[scale=1.2]{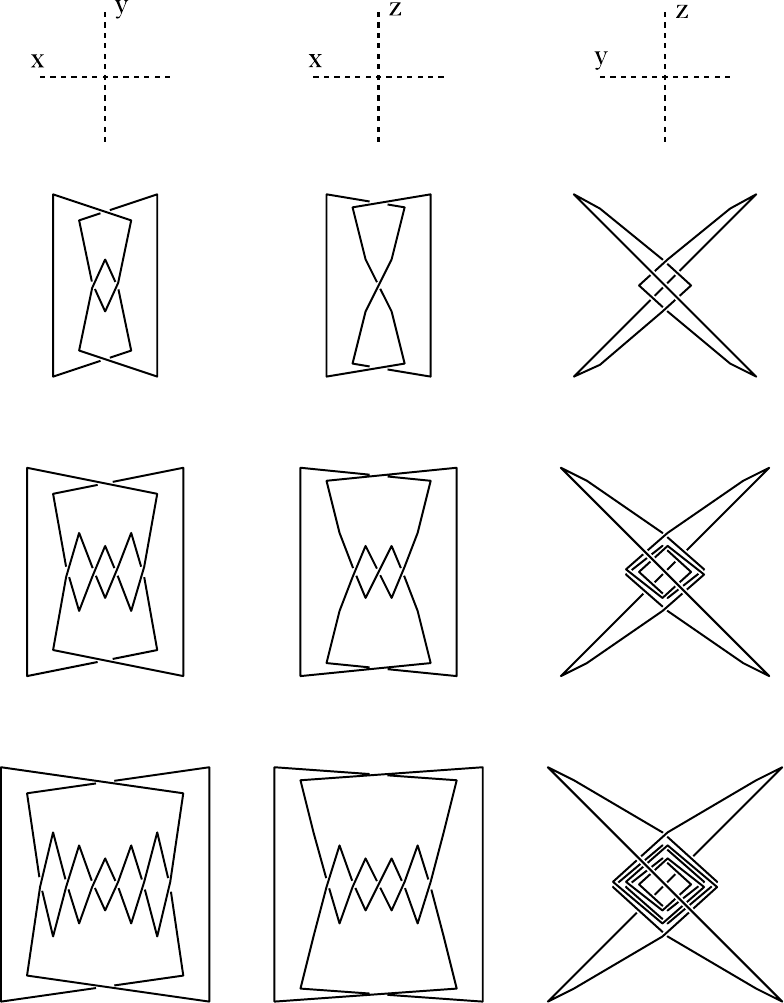}
\caption{The twist knots $C(2,-2)$, $C(4,-2)$ and $C(6,-2)$ in three orthogonal projections}
\label{projections_twist_knots_negative}
\end{figure}

\begin{remark}
(Colliding vertices) In general, vertices of polygonal knots should not project onto crossings in one of the coordinate directions. However, we make two exceptions in the twist knot family, that is, for the first and the last of the vertices in the fourth part of the polygon: 

The vertices $P_1(0,2,0)$ and $P'_1(0,-2,0)$ project to $(0,0)$ in the $x$-$z$-plane and the vertices $P_{2k+8}(-2k-2,0,0)$ and $P'_{2k+8}(2k+2,0,0)$ project to $(0,0)$ in the $y$-$z$-plane. 

This peculiarity is caused by our description of the total polygonal knot by a fourth part of it: At these points the parts are glued together. We note, that after gluing, the vertices corresponding to $P_{2k+8}$ could be omitted, because the result is a straight segment. In contrast, the vertices corresponding to $P_1$ cannot be removed, since the segments meet at an angle at these points. 
The case of the trivial knot in Figure \ref{trivial_knot} is similar and we also have vertices projecting to the central crossings.
\end{remark}

\newpage
\section{Related concepts} \label{related_concepts}
In this section we compare the simultaneous crossing number with similar concepts.

\begin{figure}[hbtp]
\centering
\includegraphics[scale=1.5]{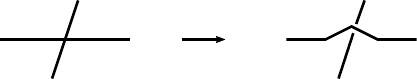}
\caption{Arbitrarily small triangle moves can be applied, so that the numbers of crossings are the same in all three orthogonal projections}
\label{triangle_move}
\end{figure}

\begin{itemize}
\item (all directions) The supercrossing index (see \cite{Denne}, Section 4.5) takes the maximal number of crossings over all directions, not just three orthogonal ones, and minimizes that value over all conformations of a given knot type. It is defined for all knots and, using projections near quadrisecants, it has been shown that for non-trivial knots it is always at least equal to 6. An averaging over all directions has also been studied (`average crossing number').
\item (no symmetry condition) Naively, one might consider the sum of crossing numbers in three orthogonal directions for \emph{all} knots, without any symmetry condition.
In this definition, the modified invariant will always have the value $3 \cdot \cn(K)$: 
				
First, construct a planar polygonal curve from a minimal knot diagram, with $\cn(K)$ crossings, and make small triangle moves at the crossings, so that the desired knot type is achieved. 
Second, rotate this construction to any plane that is not parallel to one of the coordinate planes.  This results in three knot diagrams with $\cn(K)$ crossings in each of the three orthogonal projections on the coordinate planes, see Figure \ref{triangle_move}. In this construction, the angle in the triangle move must be smaller than the smallest of the plane's angle with respect to the coordinate planes.  Therefore this naive knot invariant does not provide anything new. 
\item (two directions) The simultaneous crossing number for knots with doubly transvergent diagrams can be defined for two of the three directions. In this case, we have to decide whether (a) the two intravergent diagrams should be chosen, or (b) the doubly transvergent diagram and one of the intravergent diagrams. The latter invariant can even be defined for all strongly invertible knots.
\item (other symmetries) Besides strongly invertible knots there might be other symmetry types with three characteristic directions, allowing a definition analogous to the one discussed here. Do such symmetries exist?
\item (other invariants) Of course, other knot invariants defined by diagrams can be studied as simultaneous invariants, e.g. the simultaneous bridge number (for knots with doubly transvergent diagrams). We did not yet collect any information about them, though.
\end{itemize}

\section{A conjecture concerning the equivalence classes of strong inversions} \label{equivalence_classes}

In this section, we investigate which types of doubly transvergent diagrams occur with respect to the two equivalence classes of strong inversions.

Sakuma showed in 1986 that every invertible hyperbolic knot is strongly invertible. If, in addition, it has cyclic 
period 2, then it has two equivalence classes of strong inversions (see Sakuma \cite{Sakuma}, Proposition 3.1).
The following conjecture does not include torus knots and satellite knots. We exclude torus knots, because they have only one equivalence 
class of strong inversion, and satellite knots, because it is not clear whether all strongly invertible satellite knots with cyclic 
period 2 have doubly transvergent diagrams.

\begin{conjecture}
Let $K$ be a prime hyperbolic knot which is 
strongly invertible and has cyclic period 2.
Its two equivalence classes of strong inversions
are denoted by $\chi_1$ and $\chi_2$. 
Then, exactly one of the following conditions applies to $K$:
\begin{itemize}
\item Type I: $K$ has a doubly transvergent diagram, with one inversion axis belonging to class $\chi_1$ and the other to $\chi_2$.
\item Type II: $K$ has a doubly transvergent diagram where both axes belong to class $\chi_1$. Additionally, $K$ also has a doubly transvergent diagram where both axes belong to class $\chi_2$.
\item Type III: $K$ has a doubly transvergent diagram where both axes belong to class $\chi_1$. However, $K$ does not have a doubly transvergent diagram where both axes belong to class $\chi_2$ (or the same condition holds with the labels $\chi_1$ and $\chi_2$ switched).
\end{itemize}
\end{conjecture}

We give examples of knots for these types (taken from \cite{Sakuma}). 
Note, that the conjectured exclusiveness -- each of the example knots belongs exactly to one of the types -- is not yet proven, though.

\medskip
\noindent
\textbf{Type I}: Examples for type I are the knot $8_5$ (see Figure \ref{doubly_transvergent_diagrams}, on the left)
and two-bridge knots with fractions $p/q$ and $q^2 \not\equiv 1 (\text{mod}\;p)$.

\begin{figure}[hbtp]
\centering
\includegraphics[scale=0.8]{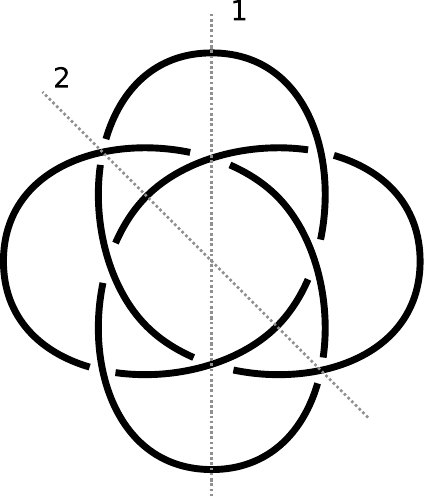}
\caption{The knot $8_{18}$. The indicated axes belong to different equivalence classes. The diagram shows period 4; 
therefore, the horizontal axis (not shown) belongs to the same equivalence class as the axis marked with 1.}
\label{diagram_8_18}
\end{figure}

\medskip
\noindent
\textbf{Type II}: An example for type II is the knot $8_{18}$. This knot has period 4 and the diagram given by Sakuma \cite{Sakuma} is doubly transvergent. 
The two orthogonal axes belong to the same equivalence class, however. An axis with different equivalence class is obtained by rotating one of the axes by $\pi/4$ (see Sakuma's diagram list, and our Figure \ref{diagram_8_18}).

\medskip
\noindent
\textbf{Type III}: Examples for type III are the (non-torus) two-bridge knots with fractions $p/q$ and $q^2 \equiv 1 (\text{mod}\;p)$,
for instance $7_4=C(4,-4)$, $7_7=C(2,2,-2,-2)$ and $9_{10}=C(4,-2,2,-4)$. See also the illustrations in Sakuma \cite{Sakuma}, Figure 3.2.

\section{Summary and outlook} \label{summary}
\subsection{Our invariant for sub-families}
We defined $\scn_{\min}(n)$ and $\scn_{\max}(n)$ as the smallest (resp. largest) simultaneous crossing number of all knots with doubly symmetric 
diagrams $K$ with $\cn(K)=n$, for $n \ge 3$. We expect, that for different classes of knots, the asymptotic behaviour of these values might
be quite different.Therefore, we give a definition with a class of knots as additional attribute.

\begin{definition} \label{subclasses_sim_min}
We denote by $\mathcal{K}_0$ the set of knots which possess doubly transvergent diagrams.
Let $\mathcal{K}_1 \subset \mathcal{K}_0$ be a subset, containing knots with arbitrarily large minimal crossing numbers.
Then we define $\scn_{\min}(\mathcal{K}_1, n)$ as the smallest simultaneous crossing number of all knots $K \in \mathcal{K}_1$ with $\cn(K)=n$.
\end{definition} 

Typical cases for $\mathcal{K}_1 \subseteq \mathcal{K}_0$ are the following: 

\begin{itemize}
	\item Composite knots,
	\item prime satellite knots,
	\item torus knots,
	\item prime hyperbolic knots with bridge number 2 (= non-torus 2-bridge knots),
	\item prime hyperbolic knots with bridge number $\ge 3$
    (possibly alternating or non-alternating).
\end{itemize}

In Definition \ref{subclasses_sim_min} we required the condition that knots with arbitrarily large minimal crossing numbers are 
contained in $\mathcal{K}_1$ because we also want to define the limit invariant:

\begin{definition} \label{subclasses_sim_limit}
If $\mathcal{K}_1 \subseteq \mathcal{K}_0$, we define 
\[
\underline{l}(\mathcal{K}_1)=\liminf_{n \to \infty}\frac{\scn_{\min}(\mathcal{K}_1, n)}{n}.
\]
\end{definition} 

This is well-defined because for the quotient we always have $\scn_{\min}(\mathcal{K}_1, n)/n \ge 3$.
For the typical cases above 
we still have to show that they contain knots with arbitrarily large minimal crossing number. 
Analogously, we may define $\scn_{\max}(\mathcal{K}_1, n)$ and $\overline{l}(\mathcal{K}_1)$.

\medskip
\subsection{A torus knot example and a conjecture}
\begin{figure}[hbtp]
\centering
\includegraphics[scale=1.8]{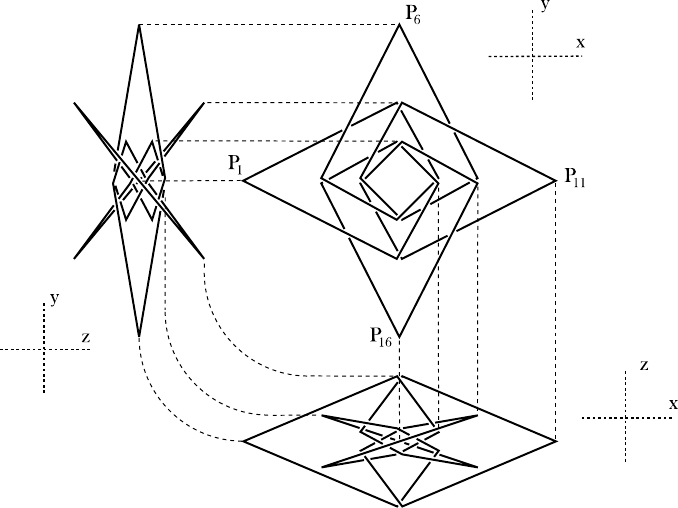}
\caption{The torus knot $T(4,5)$ shown in three orthogonal projections}
\label{projection_torus_knot_4_5}
\end{figure}

Among several examples for which we tried to minimize the total number of crossings in the three directions,
the torus knot $T(4,5)$ was the only one for which we achieved less than $4\cdot \cn(K)$ crossings.
This knot has $\cn(K) = 16$ and a minimal diagram as shown in Figure \ref{projection_torus_knot_4_5} is doubly transvergent.
The number of crossings in the projection to the x-z-plane is 19 and the number of crossings in the projection to the y-z-plane is 27.
The total number is therefore 62, and we conclude that $\scn(T(4,5)) \le 62$. (Note, that this is less than $4 \cdot 16 = 64$.)
One fourth of the polygon is described as follows: 
$P_1(-12,0,0)$, $P_2(-\delta,6,-5)$, $P_3(3,\delta,-\frac{3}{4})$, $P_4(-\delta,-3,1)$, $P_5(-6,\delta,2)$ and $P_6(0,12,0)$. 
As before, we choose $\delta=0.2$.

This leads to the following conjecture.

\begin{conjecture}
Let $\mathcal{K}_1$ be the set of torus knots which possess doubly transvergent diagrams. (This should be equal to
all torus knots $T(p,q)$ where one of the parameters is even.) Then we conjecture that $\underline{l}(\mathcal{K}_1) \le 4$.
\end{conjecture}

\medskip
\subsection{Polygonal and smooth knots}
\begin{figure}[hbtp]
\centering
\includegraphics[scale=0.15]{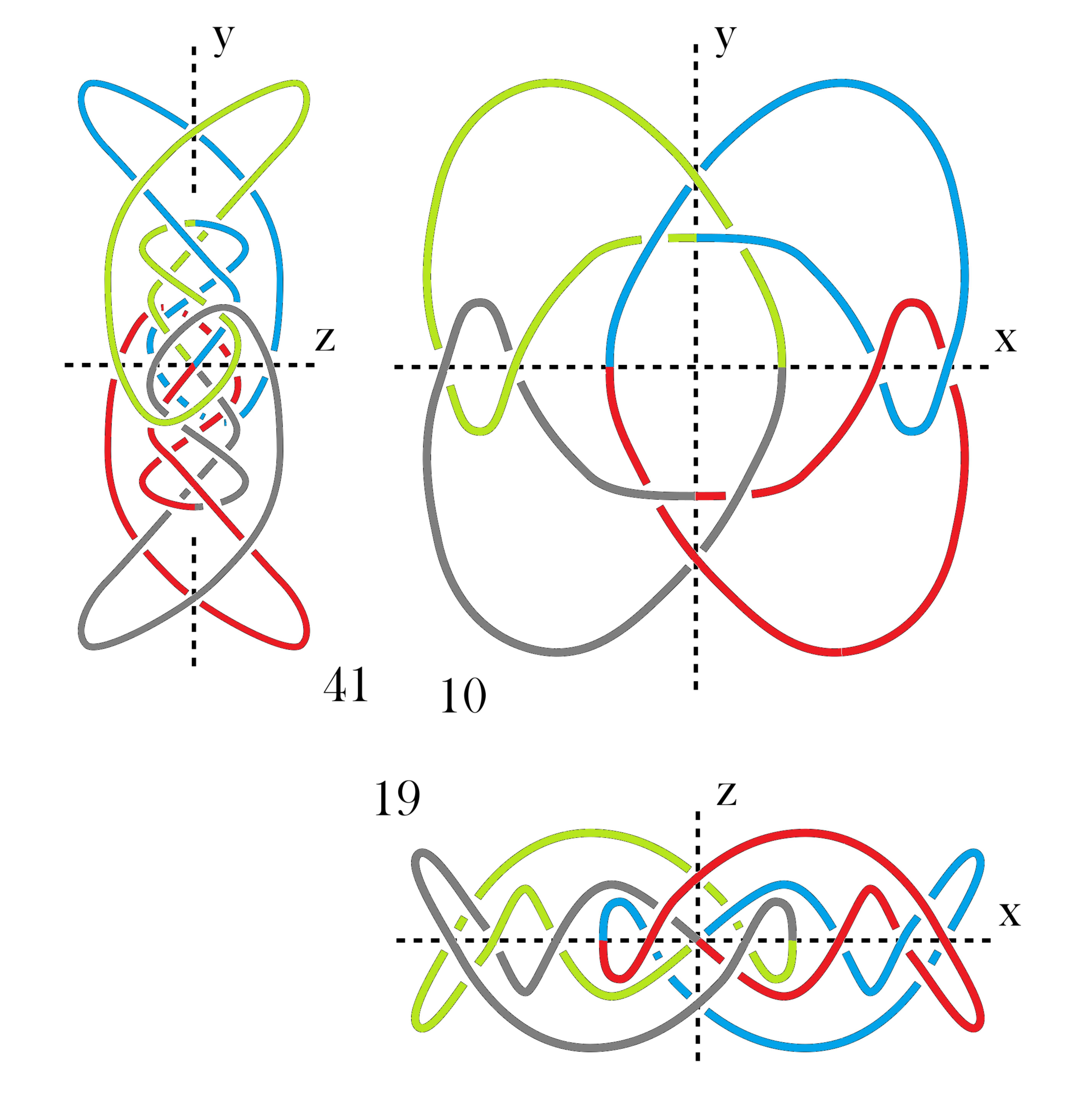}
\caption{The knot $10_{120}$ as a smooth curve, with different colours for the four parts}
\label{knot_10_120_fourColours}
\end{figure}

In Figure \ref{knot_10_120_fourColours} we show the knot $10_{120}$ in a smooth version and note, 
that without providing the vertices of a \emph{polygonal} curve, it is difficult to reconstruct the three-dimensional curve. 
Smooth curves are therefore less practicable for checking diagram properties than polygonal curves.

The three diagrams in the example have 10, 19 and 41 crossings, giving a total sum of 70 crossings.
Although this is again equal to $8 \cdot \cn(K)-10$ (see the bound in the twist knot case in Theorem \ref{twist_knot_theorem}),
for other alternating knots we also found upper bounds which are smaller than $8 \cdot \cn(K)-10$, e.g. $\scn(8_{18})\le 50$, $\scn(9_{35})\le 52$ and $\scn(10_{136})\le 54$.

Currently, we do not have a conjecture for the value of $\underline{l}$ for alternating prime hyperbolic knots with bridge number $\ge 3$. For non-torus two-bridge knots, the bound in Theorem \ref{twist_knot_theorem} could be optimal and then we would have $\underline{l}=8$ for prime hyperbolic knots with bridge number 2.

\medskip
\subsection{Random knots} The following project might be interesting for researchers who already have experience with generating and identifying random knots.

\begin{project}
Generate randomly fourth parts of polygonal knots with doubly symmetric diagrams. In each case, construct the complete knot and determine the numbers of crossings in the three directions. If also the identification of the knot type is successful (using Knotscape or SnapPy, resulting in a knot with 16 or less crossings), then add the knot type and its polygonal description to a result list, which can be searched for interesting examples.
\end{project}

\clearpage
\newpage

\section{Appendix A, on the symmetry of chord diagrams} \label{appendix_a}
We mentioned that the question whether $\scn(3_1)=14$ could be decided by an enumeration of all cases with a total number of crossings less than $14$. A possible tool for that is the analysis of the symmetries of the chord diagrams, taken in the three directions simultaneously.

\begin{figure}[hbtp]
\centering
\includegraphics[scale=0.95]{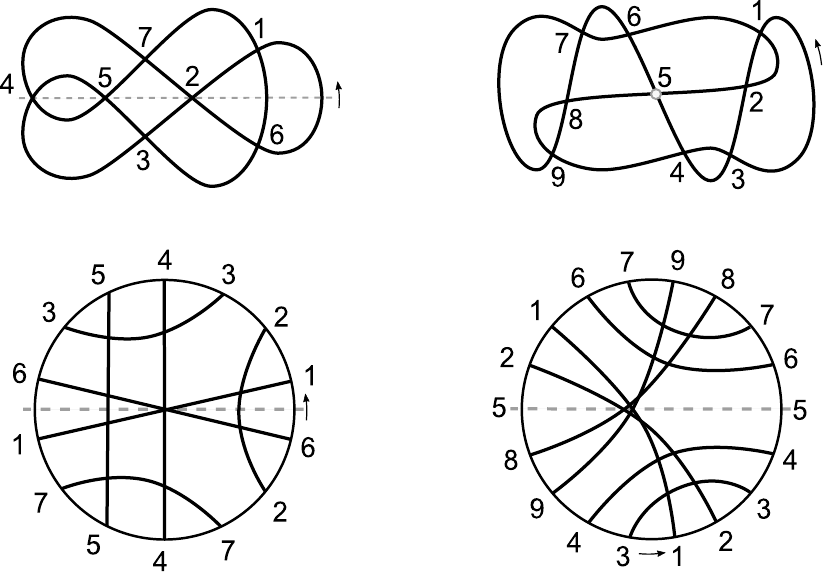}
\caption{Shadows of a transvergent diagram for the knot $7_6$ and of an intravergent diagram for the knot $7_7$, together with their chord diagrams.}
\label{chord_diagrams}
\end{figure}

We start with two examples of chord diagrams. One for a transvergent diagram, the other for an intravergent diagram, see Figure \ref{chord_diagrams}.

The axes of rotation are shown in gray. The effect on the chord diagrams is a mirror symmetry with regard to these axes. In the intravergent case, one of the chords coincides with the rotation axis.

\begin{figure}[hbtp]
\centering
\includegraphics[scale=1]{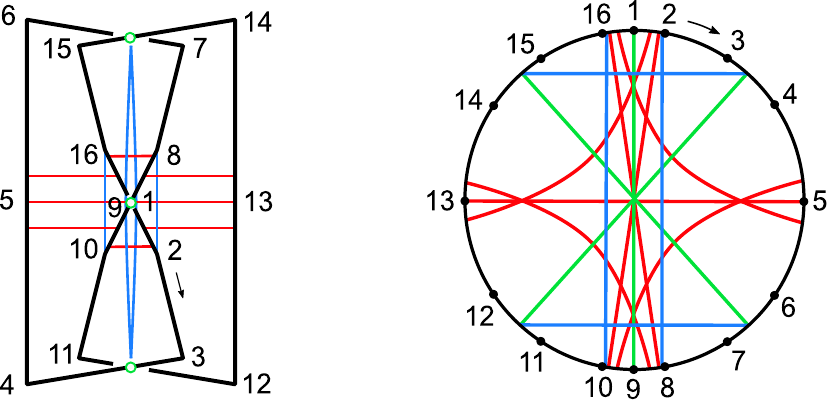}
\caption{Simultaneous chord diagrams for the trefoil: Blue chords correspond to the doubly transvergent diagram, red and green chords to the intravergent diagrams.}
\label{chord_diagrams_trefoil}
\end{figure}

The chord diagram of a doubly transvergent diagram has two orthogonal symmetry axes. For the trefoil, given as $C(2,-2)$ as in our parametrization in Section \ref{asymptotic_behaviour} (with details in Appendix C), the three chord diagrams can be seen simultaneously in Figure \ref{chord_diagrams_trefoil}. 

In this chord diagram, the numbers refer to the vertices of the polygon, not to the crossings in the diagrams. Note, that the blue and green chords meeting between vertices 11 and 12 have the exact same end points on the polygon. On the other hand, the blue and red chords meeting close to vertex 10 do not have coinciding end points.
It is not yet clear, how chord diagrams can help us in the enumeration problem. 

\clearpage
\newpage

\section{Appendix B, workflow for using the `Blender' software} \label{appendix_b}
Optimizing a three-dimensional curve to achieve a small total crossing number requires a 3D software tool. We explain for beginners how to use the `Blender' software.

\small
\medskip
\noindent
Preparation and two-dimensional modelling:
\begin{itemize}
    \item Install and open Blender. You see the standard scene containing a cube. We don't need this cube, so delete it.
    \item In the top left, use the `Add' menu, choose `curve' and `Bézier' to   create a Bézier curve.
    \item In the top right use the coloured coordinate symbol and press on $z$, to view the scene from above. Use the mouse wheel to get nearer. 
    \item Object and edit modes: The curve is shown in object mode. Go to the button in the top left and change it to `Edit mode'. The curve now has Bézier handles for the two vertices.
    \item Clicking on one of the vertices activates it for operating: moving (press g/G), rotating (press R) and scaling (press S). The operation is done with the mouse and finished with left clicking it. When you finish with a right click, the operation is not processed.
    \item Elongate your curve by pressing E. This creates a new connected vertex.
    \item Now model a two-dimensional doubly transvergent knot diagram. You can close the curve by activating the last and the first vertex (shift-click) and then press F.
    \item You can assign a radius to this curve as follows: Find the curve symbol in the lower right side. Go to sub-menu `Geometry' and `Bevel' and increase the `Depth' attribute.
\end{itemize}

\noindent
Three-dimensional modelling and optimization:
\begin{itemize}
    \item Try moves and rotations in the modelling space by moving the mouse while pressing the mouse wheel, and do the same while also pressing shift.
    \item A useful feature are simultaneous windows for different views of the curve: In the upper left corner you can open a sub-window (cursor changes to a cross symbol). After a second window is opened you can modify the view direction by pressing e.g. $x$ in the coordinate symbol.
    \item Blender has a function for symmetric constructions, but not for our rotationally symmetric case. You therefore have to make sure that the curve has the required symmetries yourselves. (It would be possible to program an add-on which constructs the total curve when only one fourth is drawn.)
    \item Your knot is still two-dimensional and you have to move the vertices in a way that the desired over- and under-crossings are achieved. To do that, you probably need more vertices to be more flexible. This can be done by choosing `Subdivide' from the `Segments' menu.
    \item Do the three-dimensional modelling while respecting the symmetry requirements. This results in a knot diagram looking fine from the top but probably terrible from the other two orthogonal directions. It is necessary to simplify these views while leaving the knot type unchanged. At the same time you can try to optimize the numbers of crossings.
    \item Moving vertices up or down can be easily done without the mouse: Activate the vertex, then press G, z, 1, for moving it up by one unit. A symmetric vertex is moved down with G, z, -1.
\end{itemize}

\noindent
Clean-up and obtaining final polygonal data:
\begin{itemize}
    \item In Blender, often the curve is transformed to a mesh object. This is not necessary for our purpose.
    \item To transform the smooth curve into a polygon, you press A in order to select all vertices, then use the right mouse button and choose `Set Spline Type', `Poly'.
    \item To be able to check the symmetry rigorously, you should try to move the vertices to integer coordinates. The coordinates for one fourth of the curve can be used for an input script in Python, similar to our pseudo code examples. If this input creates your knot, then you are sure that the coordinates are correct.
    \item We don't need materials, lighting and rendering for these modelling activities.
\end{itemize}

\normalsize
More advanced material can be found at the web site of Robert Lipshitz (University of Oregon). 
He mentions Blender tutorials, and explains the use of line art and Python scripting.

\section{Appendix C, on the remaining twist knot cases} \label{appendix_c}
Section \ref{asymptotic_behaviour} contained the pseudo code for generating (one fourth of) the polygons for $C(2k,2)$ for odd $k$.
The following algorithm gives the coordinates for $C(2k,2)$ for even $k$.

\footnotesize
\medskip
\begin{verbatim}
Case 2: C(2k,2) for even k

input m    # Examples: m=1 and k=2 give C(4,2)
k = 2*m    #           m=2 and k=4 give C(8,2)
d = 0.2    # The parameter delta
  
P[1] = (0, 2, 0)
P[2] = (1,-d,-2)
P[3] = (2,-3, d)
P[4] = (3, d, 3)

j = 4      # counter for vertices, starting with P[5]

for i = 0 to m-2:  # spiral part enlargement
    x = 4*i+4; y = 2*i+4; z = -d; j = j+1; P[j] = (x,y,z)    		
    x = x+1;   z =  -y;   y = -d; j = j+1; P[j] = (x,y,z)		
    x = x+1;   y = z-1;   z =  d; j = j+1; P[j] = (x,y,z)		
    x = x+1;   z =  -y;   y =  d; j = j+1; P[j] = (x,y,z)
end

j = j+1; P[j] = ( 2*k  ,k+2,  -d)
j = j+1; P[j] = ( 2*k+1,k+5,-k-4)
j = j+1; P[j] = (-2*k-2,k+6,-k-6)
j = j+1; P[j] = (-2*k-2,  0,   0)
\end{verbatim}

\normalsize
\medskip
Note, that the coordinates with index $\le 4$ are nearly identical to the spiral part block evaluated for $i=-1$.
The only difference is, that we have $P[1]=(0,2,0)$ and not $P[1]=(0,2,-\delta)$, with the effect that this point
is projected to the origin in the projection to the $x$-$z$-plane. This is necessary to satisfy the symmetry condition.

\normalsize
\medskip
The following algorithm in pseudo code gives the coordinates for $C(2k,-2)$ for odd $k$:

\footnotesize
\medskip
\begin{verbatim}
Case 3: C(2k,-2) for odd k

input m    # Examples: m=0 and k=1 give C(2,-2)
k = 2*m+1  #           m=1 and k=3 give C(6,-2)
d = 0.2    # The parameter delta
  
P[1] = (0, 2, 0)
P[2] = (1,-d,-2)

j = 2      # counter for vertices, starting with P[3]

for i = 0 to m-1:  # spiral part enlargement
    x = 4*i+2; y = -(2*i+3); z =  d; j = j+1; P[j] = (x,y,z)    		
    x = x+1;   z =  -y;      y =  d; j = j+1; P[j] = (x,y,z)		
    x = x+1;   y = z+1;      z = -d; j = j+1; P[j] = (x,y,z)		
    x = x+1;   z =  -y;      y = -d; j = j+1; P[j] = (x,y,z)
end

j = j+1; P[j] = ( 2*k  ,-k-4,-k-5)
j = j+1; P[j] = (-2*k-2,-k-6,-k-6)
j = j+1; P[j] = (-2*k-2,   0,   0)
\end{verbatim}

\bigskip
\normalsize
The following algorithm gives the coordinates for $C(2k,-2)$ for even $k$.
It combines properties of cases 2 and 3.

\footnotesize
\medskip
\begin{verbatim}
Case 4: C(2k,-2) for even k

input m    # Examples: m=1 and k=2 give C(4,-2)
k = 2*m    #           m=2 and k=4 give C(8,-2)
d = 0.2    # The parameter delta
  
P[1] = (0, 2, 0)
P[2] = (1,-d,-2)
P[3] = (2,-3, d)
P[4] = (3, d, 3)

j = 4      # counter for vertices, starting with P[5]

for i = 0 to m-2:  # spiral part enlargement
    x = 4*i+4; y = 2*i+4; z = -d; j = j+1; P[j] = (x,y,z)    		
    x = x+1;   z =  -y;   y = -d; j = j+1; P[j] = (x,y,z)		
    x = x+1;   y = z-1;   z =  d; j = j+1; P[j] = (x,y,z)		
    x = x+1;   z =  -y;   y =  d; j = j+1; P[j] = (x,y,z)
end

j = j+1; P[j] = ( 2*k  ,k+4,k+5)
j = j+1; P[j] = (-2*k-2,k+6,k+6)
j = j+1; P[j] = (-2*k-2,  0,  0)
\end{verbatim}

\begin{figure}[hbtp]
\centering
\includegraphics[scale=0.22]{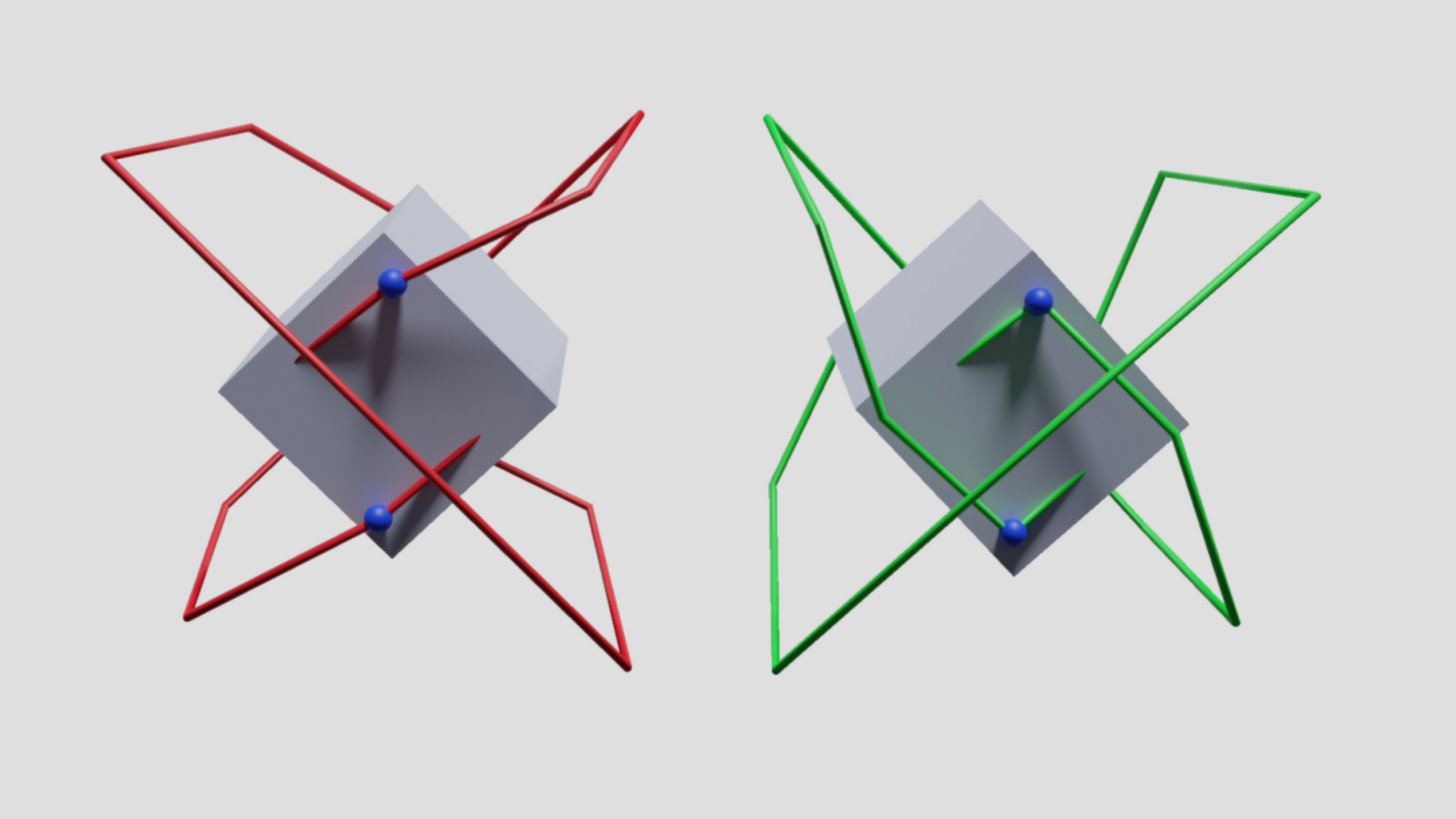}
\caption{The polygons outside of the spiral part are shown for the family $C(2k,-2)$ on the left, and for $C(2k,2)$ on the right}
\label{boxed_center}
\end{figure}

\normalsize
\bigskip
In Figure \ref{boxed_center} we illustrate that the central spiral parts are identical for the cases $C(2k,2)$ and $C(2k,-2)$.
By masking the identical spiral parts in a cube, we emphasize the differences in the outer polygonal parts.
The vertices marked in blue separate the outer polygonal parts, which are different, from the central parts.

The fact that the spiral parts are identical in all four cases is easier to see in the following unified version of the algorithm. We write $C(2k,2\rho)$ for $C(2k,2)$ and $C(2k,-2)$, with $\rho = \pm 1$ (see lines 3 and 4), and introduce a parameter $e = \pm 1$ (line 13). This parameter allows a simpler description of the spiral part. The vertices marked in blue in Figure \ref{boxed_center} are the last ones in the while loop; only one case distinction remains (lines 22--27). Note, that $\rho$ affects the $z$-coordinate in line 30.

\footnotesize
\medskip
\begin{verbatim}
 1  Unified code: C(2k,2*rho)
 2 
 3  input k    # k=3, rho=1  give C(6,2)
 4  input rho  # k=3, rho=-1 give C(6,-2)
 5  d = 0.2    # The parameter delta
 6   
 7  P[1] = (0, 2, 0)
 8  P[2] = (1,-d,-2)
 9
10  # The second vertex is the starting point of the while loop
11  x = 1; y = -d; z = -2
12 
13  e = -1     # negative values indicate decreasing values in y or z direction
14  j =  2     # counter for vertices, starting with P[3]
15 
16  while x < 2*k-1:  # spiral part enlargement
17     e = -e
18     x = x+1; y = z-e; z = e*d; j = j+1; P[j] = (x,y,z)    		
19     x = x+1; z =  -y; y = e*d; j = j+1; P[j] = (x,y,z)		
20  end
21 
22  if rho==1:  # case C(2k,2)
23     j = j+1; P[j] = (2*k  ,e*(k+2),-e*d)
24     j = j+1; P[j] = (2*k+1,e*(k+5),-e*(k+4))
25  
26  if rho==-1: # case C(2k,-2)
27     j = j+1; P[j] = (2*k  ,e*(k+4), e*(k+5))
28 
29  # no case distinction for the final two vertices
30  j = j+1; P[j] = (-2*k-2,e*(k+6),-e*rho*(k+6))
31  j = j+1; P[j] = (-2*k-2,      0,           0)
\end{verbatim}

\normalsize

\vspace{2cm}

\end{document}